\documentclass[12pt]{amsart}
\usepackage{t1enc}
\usepackage[latin2]{inputenc}
\usepackage{amsmath}
\usepackage{amsthm}
\usepackage{amssymb}
\usepackage{stmaryrd}
\usepackage[T1]{fontenc}
\usepackage{enumerate}
\usepackage{verbatim}
\usepackage{amscd}

\newif\ifcommented
\commentedtrue

\newcommand{\comm}[1]{}

\ifcommented
\renewcommand{\comm}[1]{\fbox{\fbox{\begin{minipage}{300pt}#1\end{minipage}}}}

\fi

\newtheorem{Thm}{Theorem}[section]
\newtheorem{Prop}[Thm]{Proposition}

\newtheorem{Lem}[Thm]{Lemma}

\newtheorem{Kov}[Thm]{Corollary}

\newtheorem{Prob}[Thm]{Problem}

\newtheorem{Fact}[Thm]{Fact}

\newtheorem{claim}{Claim}[Thm]

\newtheorem{obs}[Thm]{Observation}

\theoremstyle{definition}

\newtheorem{definition}[Thm]{Definition}

\newcommand{\mc}{\mathcal}
\newcommand{\mf}{\protect\mathfrak}
\newcommand{\mbb}{\mathbb}

\newcommand{\PP}{\mathbb{P}}

\newcommand{\ra}{\rangle}
\newcommand{\la}{\langle}

\newcommand{\0}{\emptyset}
\newcommand{\bs}{\backslash}
\newcommand{\vd}{\Vdash}
\newcommand{\Ra}{\Rightarrow}
\newcommand{\fv}{\rightarrow}

\newcommand{\al}{\alpha}
\newcommand{\be}{\beta}

\newcommand{\de}{\delta}
\newcommand{\eps}{\varepsilon}
\newcommand{\lam}{\lambda}
\newcommand{\ka}{\kappa}

\newcommand{\om}{\omega}
\newcommand{\bc}{\begin{center}}
\newcommand{\ec}{\end{center}}
\newcommand{\fel}{\ge}

\newcommand{\dom}{\mathop{\mathrm{dom}}\nolimits}
\newcommand{\ran}{\mathop{\mathrm{ran}}\nolimits}

\newcommand{\zfc}{\mathop{\mathsf{ZFC}}\nolimits}

\newcommand{\ch}{\mathop{\mathsf{CH}}\nolimits}

\newcommand{\acal}{{\mathcal{A}}}
\newcommand{\ical}{{\mathcal{I}}}
\newcommand{\jcal}{{\mathcal{J}}}

\newcommand{\pcal}{{\mathcal{P}}}

\newcommand{\zcal}{{\mathcal{Z}}}

\makeindex

\newif\ifdeveloping

\developingtrue

\def\myheads#1;#2;{
\pagestyle{myheadings}
\markboth{{\sc\hfill #1\hfill\protect\makebox[0cm][r]{\rm\today}}}
{{\sc\protect\makebox[0cm][l]{\rm\today}\hfill #2\hfill}}
}

\def\<{\left\langle}
\newcommand{\muv}{{\vec \mu}}
\def\>{\right\rangle}

\newcommand{\subs}{\subseteq}
\newcommand{\setm}{\setminus}

\def\br#1;#2;{\bigl[ {#1} \bigr]^ {#2} }
\newcommand{\fini}{\mathsf{fin}}
\newcommand{\lerk}{\le_{\rm RK}}
\newcommand{\lerb}{\le_{\rm RB}}

\begin{document}

\title{More on cardinal invariants of analytic P-ideals}
\thanks{The preparation of this paper was supported by
Hungarian National Foundation for Scientific Research grant no 61600, 68262
and 63066}

\author[B. Farkas]{Barnab\'as Farkas}
\address{Budapest University of Technology and Economics (BME)}
\email{barnabasfarkas@gmail.com}


\author[L. Soukup]{Lajos Soukup}
\address{Alfr{\'e}d R{\'e}nyi Institute of Mathematics }
\email{soukup@renyi.hu}

\begin{abstract}
Given an ideal $\mc{I}$ on $\om$ let
$\mf{a}(\mc{I}) $  ($\bar{\mf{a}}(\mc{I})$)
be minimum of the cardinalities of infinite (uncountable)
maximal $\ical$-almost disjoint subsets of $[{\omega}]^{\omega}$.
We show that $\mf{a}({\ical_h})>{\omega}$
if $\ical_h$ is a summable ideal; but
$\mf{a}({\zcal_{\muv}})={\omega}$ for any tall density ideal
$\zcal_{\muv}$
including the density zero ideal $\zcal$.
On the other hand,  you have $\mf{b}\le \bar{\mf{a}}({\ical})$ for any analytic
$P$-ideal $\ical$, and
$\bar{\mf{a}}({\zcal_{\muv}})\le \mf{a}$
for each density ideal $\zcal_\muv$.

For each ideal $\ical$ on $\om$ denote $\mf{b}_\ical$ and $\mf{d}_\ical$ the unbounding and dominating numbers of $\la\om^\om,\le_\ical\ra$ where $f\le_\ical g$ iff $\{n\in\om:f(n)>g(n)\}\in\ical$. We show that $\mf{b}_\ical=\mf{b}$ and $\mf{d}_\ical=\mf{d}$ for each analytic P-ideal $\ical$.

Given a Borel ideal $\mc{I}$ on ${\omega}$
we say that a poset $\PP$ is {\em $\mc{I}$-bounding }
iff $\forall x\in \mc{I}\cap V^{\PP}$
$\exists y\in \mc{I}\cap V$  $x\subseteq y$.
$\PP$ is {\em $\mc{I}$-dominating}  iff
$\exists y\in \mc{I}\cap V^{\PP}$
$\forall x\in \mc{I}\cap V$ $x\subseteq^* y$.

For each analytic P-ideal $\mc{I}$
if a poset $\PP$ has the Sacks pro\-per\-ty then $\PP$ is
$\ical$-bounding; moreover if $\ical$ is tall as well then
the property $\ical$-bounding/$\ical$-dominating implies
${\omega}^{\omega}$-bounding/adding dominating reals, and the converses of these two implications are false.

For the density zero ideal $\zcal$ we can prove more:
(i) a poset $\PP$ is $\zcal$-bounding iff it has the
Sacks property,
(ii) if $\PP$ adds a slalom capturing all ground model reals then
$\PP$ is $\mc{Z}$-dominating.

\end{abstract}

\maketitle

\section{Introduction}

In this paper we
  investigate some properties of some cardinal
invariants associated with  analytic $P$-ideals. Moreover we
   analyze  related ``bounding'' and ``dominating'' properties of forcing
  notions.

Let us denote $\mathsf{fin}$ the Frechet ideal on ${\omega}$, i.e.
$\fini=[\om]^{<\om}$.
Further we always assume that if $\ical$ is an ideal on $\om$ then
the ideal is {\em proper}, i.e. $\om\notin \ical$, and $\fini\subseteq\ical$,
so especially $\mc{I}$ is {\em non-principal}. Write $\ical^+=\mc{P}(\om)\bs\ical$ and $\ical^*=\{\om\bs X: X\in\ical\}$.

An ideal $\ical$ on $\om$ is {\em analytic} if $\ical\subseteq\mc{P}(\om)\simeq 2^\om$ is an analytic set in the usual product topology. $\ical$ is a {\em P-ideal}
if for each countable $\mc{C}\subseteq\mc{I}$ there is an $X\in\mc{I}$ such that $Y\subseteq^* X$ for each $Y\in\mc{C}$, where $A\subseteq^* B$ iff $A\bs B$ is finite. $\ical$ is {\em tall} (or {\em dense}) if each infinite subset of $\om$ contains an infinite element of $\ical$.

A function $\varphi:\pcal({\omega})\to [0,\infty]$
is a {\em submeasure on ${\omega}$} iff
${\varphi}(X)\le \varphi(Y)$ for $X\subs Y\subs {\omega}$,
$\varphi(X\cup Y)\le \varphi(X)+\varphi(Y)$ for $X,Y\subs {\omega}$, and
$\varphi(\{n\})<\infty$ for $n\in {\omega}$.
A submeasure $\varphi$ is {\em lower semicontinuous}
iff $\varphi(X)=\lim_{n\fv\infty}\varphi(X\cap n)$ for each $X\subseteq\om$. A submeasure $\varphi$ is {\em finite} if $\varphi(\om)<\infty$. Note that if $\varphi$ is a lower semicontinuous submeasure on $\om$ then $\varphi(\bigcup_{n\in\om}A_n)\le\sum_{n\in\om}\varphi(A_n)$ holds as well for $A_n\subseteq\om$.
We assign the {\em exhaustive  ideal $\mathrm{Exh}(\varphi)$} to a submeasure $\varphi$ as follows
\begin{displaymath}
\mathrm{Exh}(\varphi)=\big\{X\subseteq {\omega}:\lim_{n\fv\infty}\varphi(X\bs n)=0\big\}.
\end{displaymath}

Solecki, \cite[Theorem 3.1]{So}, proved  that an ideal $\ical\subs \pcal({\omega})$
is an analytic $P$-ideal or $\ical=\mc{P}(\om)$ iff $\ical=\mathrm{Exh}(\varphi)$ for some
lower semicontinuous finite submeasure.
Therefore  each analytic P-ideal is $F_{\sigma\delta}$ (i.e. $\Pi^0_3$) so a Borel subset of $2^\om$.
It is straightforward to see that if $\varphi$ is a lower semicontinuous finite submeasure on $\om$ then the ideal $\mathrm{Exh}(\varphi)$ is tall iff $\lim_{n\fv\infty}\varphi(\{n\})=0$.

Let $\ical$ be an ideal on $\om$. A family $\acal\subseteq\ical^+$ is {\em $\ical$-almost-disjoint}
($\ical$-AD in short),   if $A\cap B\in\ical$ for each $\{A,B\}\in
[\acal]^2$. An  $\ical$-AD family $\acal$ is an
{\em $\ical$-MAD family} if for each $X\in\ical^+$
there exists an $A\in\acal$
such that $X\cap A\in\ical^+$, i.e.
$\acal$ is $\subseteq$-maximal among the $\ical$-AD families.

Denote $\mf{a}(\ical)$ the
minimum of the cardinalities of infinite $\ical$-MAD families.
In Theorem \ref{tm:a} we show that $\mf{a}({\ical_h})>{\omega}$
if $\ical_h$ is a summable ideal; but
$\mf{a}({\zcal_{\muv}})={\omega}$ for any tall density ideal
$\zcal_{\muv}$
including the {\em density zero ideal}
\[ \zcal=\Big\{A\subseteq\om:\lim_{n\fv\infty}\frac{|A\cap n|}{n}=0\Big\}.\]
On the other hand, if you define
 $\bar{\mf{a}}(\ical)$ as
minimum of the cardinalities of uncountable $\ical$-MAD families
then you have $\mf{b}\le \bar{\mf{a}}({\ical})$ for any analytic
$P$-ideal $\ical$, and
$\bar{\mf{a}}({\zcal_{\muv}})\le \mf{a}$
for each density ideal $\zcal_\muv$
(see Theorems \ref{a0a} and \ref{ba0}).

In Theorem \ref{cohen_indest} we prove under $\ch$ the existence of
an uncountable  Cohen-indestructible $\mc{I}$-MAD families for each analytic P-ideal
$\mc{I}$.

A sequence $\la A_\al:\al<\ka\ra\subset[\om]^\om$ is a {\em tower} if
it is
{\em $\subseteq^*$-descending,} i.e. $A_\be\subseteq^* A_\al$ if
$\al\le\be<\ka$, and it has no {\em pseudointersection}, i.e. a set
$X\in[\om]^\om$ such that $X\subseteq^* A_\al$ for each $\al<\ka$.
In Section \ref{sec:tower} we show it is consistent that the continuum
is arbitrarily large and
for each  tall analytic $P$-ideal $\mc{I}$ there is  towers of
height $\om_1$ whose elements are in  $\mc{I}^*$.

Given an ideal $\ical$ on $\om$ if $f,g\in\om^\om$ write $f\le_\ical
g$ if $\{n\in\om:f(n)>g(n)\}\in\ical$.
As usual let $\le^*=\le_{\mathsf{fin}}$.
The unbounding and dominating numbers of the partially ordered set
$\la\om^\om,\le_\ical\ra$, denoted by $\mf{b}_\mc{I}$ and
$\mf{d}_\mc{I}$ are defined in the natural way, i.e.
$\mf{b}_\mc{I}$ is the minimal size of a $\le_{\mc{I}}$-unbounded family, and
$\mf{d}_\mc{I}$ is the minimal size of a $\le_{\mc{I}}$-dominating family.
By these notations $\mf{b}=\mf{b}_{\mathsf{fin}}$ and
$\mf{d}=\mf{d}_{\mathsf{fin}}$.
In Section \ref{sec:bd} we show that $\mf{b}_\mc{I}=\mf{b}$ and
$\mf{d}_\mc{I}=\mf{d}$ for each analytic P-ideal $\ical$. We also
prove, in Corollary \ref{cor:dom}, that for any analytic P-ideal $\mc{I}$
 a poset $\PP$  is $\le_{\mc{I}}$-bounding  iff it is
 ${\omega}^{\omega}$-bounding, and $\PP$  adds
 $\le_{\mc{I}}$-dominating reals iff it adds
dominating reals.

In Section \ref{sec:forcing} we introduce the $\mc{I}$-bounding and
$\mc{I}$-dominating properties of forcing notions for Borel ideals:
$\PP$ is {\em $\mc{I}$-bounding} iff any  element of
$\mc{I}\cap V^{\PP}$ is contained
in some element of $\mc{I}\cap V$;
$\PP$ is {\em $\mc{I}$-dominating} iff there is  an  element in
$\mc{I}\cap V^{\PP}$ which mod-finite contains
all elements of $\mc{I}\cap V$.

In Theorem \ref{ibounddom}
we show that for each tall analytic P-ideal $\mc{I}$ if
a forcing notion is $\mc{I}$-bounding then it is $\om^\om$-bounding,
and if it is $\mc{I}$-dominating then it adds dominating reals.
Since the random real forcing is not $\ical$-bounding for each tall summable and tall density ideal $\ical$ by Proposition \ref{random}, the converse of the first implication is false.
Since a ${\sigma}$-centered forcing can not be $\mc{I}$-dominating
for a tall analytic P-ideal $\mc{I}$
by Theorem \ref{tm:sigma},
the standard dominating real forcing $\mathbb D$ witnesses that the
converse of the second implication is also false.

We prove in Theorem \ref{sacksibound} that the Sacks property implies the $\mc{I}$-bounding
property for each analytic P-ideal $\mc{I}$.

Finally, based on a theorem of Fremlin we show  that the
$\mc{Z}$-bounding property is equivalent to the Sacks property.

\section{Around the almost disjointness number of ideals}\label{sec:a}

For any ideal $\ical$ on ${\omega}$
denote $\mf{a}(\ical)$ the
minimum of the cardinalities of infinite $\ical$-MAD families.

To start the investigation of this cardinal invariant
 we recall the definition of two special classes of
analytic $P$-ideals: the density ideals and the summable ideals (see \cite{Fa}).

\begin{definition}
 Let $h:\om\fv\mbb{R}^+$ be a function such that $\sum_{n\in\om}
 h(n)=\infty$. The {\em summable ideal corresponding to $h$} is
\[ \ical_h=\Big\{A\subseteq\om:\sum_{n\in A} h(n)<\infty\Big\}.\]

Let $\<P_n:n<{\omega}\>$ be a decomposition of $\om$ into pairwise
disjoint nonempty finite sets and let $\muv=\<\mu_n:n\in\om\>$ be a
sequences of probability measures, $\mu_n:\mc{P}(P_n)\fv[0,1]$. The {\em density ideal generated by $\muv$}
is
\[ \zcal_{\muv}=\big\{A\subseteq\om:\lim_{n\fv\infty}\mu_n(A\cap P_n)=0\big\}.\]
\end{definition}

A summable ideal $\ical_h$ is tall iff
$\lim_{n\fv\infty}h(n)=0$; and a density ideal $\zcal_{\vec{\mu}}$ is
tall iff
\[\tag{\dag} \lim_{n\fv\infty}\max_{i\in P_n}\mu_n(\{i\})=0.\]

Clearly the density zero ideal $\zcal$ is a tall density ideal, 
and the summable and
the density ideals are proper ideals.

\begin{Thm}\label{tm:a}
  \begin{enumerate}[(1)]
   \item  $\mf{a}({\ical_h})>{\omega}$ for any summable ideal $\ical_h$.
   \item   $\mf{a}({\zcal_{\muv}})={\omega}$ for any tall density ideal $\zcal_\muv$.
  \end{enumerate}
\end{Thm}

\begin{proof}
(1):
We show that if $\{A_n:n<{\omega}\}\subseteq \ical^+_h$ is
$\ical$-AD then there is  $B\in \ical^+_h$ such that
$B\cap A_n\in\ical$ for $n\in {\omega}$.

For each $n\in {\omega}$ let $B_n\subseteq
A_n\setminus\cup\{A_m:m<n\}$ 
be finite such that
$\sum_{i\in B_n}h(i)>1$, and put
\begin{displaymath}
B=\cup\{B_n:n\in {\omega}\}.
\end{displaymath}

(2): Write $\muv=\la\mu_n:n\in\om\ra$ and $\mu_n$ concentrates on
$P_n$. By $(\dag)$ we have $\lim_{n\fv\infty}|P_n|=\infty$.

Now  for each $n$ we can choose $k_n\in {\omega}$ and a partition
$\{P_{n,k}:k<k_n\}$ of $P_n$ such that
  \begin{itemize}
  \item[(a)] $\lim_{n\fv\infty} k_n=\infty$,
\item[(b)] if $k<k_n$ then ${\mu}_n(P_{n,k})\ge \frac {1}{2^{k+1}}$.

  \end{itemize}
Put $A_k=\cup\{P_{n,k}:k<k_n\}$ for each $k\in\om$. We show that
$\{A_k:k\in\om\}$ is a $\zcal_\muv$-MAD family.

If $k_n>k$ then ${\mu}_n(A_k\cap P_n)={\mu}_n(P_{n,k})\ge
\frac {1}{2^{k+1}}$. 
Since for an arbitrary $k$ for all but finitely many $n$ we have
$k_n>k$ it follows that 
\[\limsup_{n\fv\infty} {\mu}_n(A_k\cap P_n)=
\limsup_{n\fv\infty}{\mu}_n(P_{n,k})\ge \limsup_{n\fv\infty}\frac
{1}{2^{k+1}}= \frac{1}{2^{k+1}}>0,\]
thus $A_k\in \zcal_{\muv}^+$.

Assume that $X\in \zcal_{\muv}^+$. Pick   $\varepsilon>0$ with 
$\limsup_{n\fv\infty} {\mu}_n(X\cap P_n)>\varepsilon$. For a large
enough $k$ we have $\frac{1}{2^{k+1}}<\frac{\varepsilon}{2}$ so if
$k<k_n$ then
\[ {\mu}_n(P_n\setminus \cup\{P_{n,i}:i\le k\})\le \frac{1}{2^{k+1}}<\frac{\varepsilon}{2}.\]
So for each large enough $n$ there is $i_n\le k$ such that
${\mu}_n(X\cap P_{n,i_n})>\frac{{\varepsilon}}{2(k+1)}$. Then
$i_n=i$ for infinitely many $n$, so $\limsup_{n\fv\infty}
{\mu}_n(X\cap A_i)\ge \frac{{\varepsilon}}{2(k+1)}$, and so
$X\cap A_i\in\zcal^+_\muv$.
\end{proof}

This Theorem gives  new proof of the following well-known fact:
\begin{Kov}
The density zero ideal $\zcal$ is not  a summable ideal.
\end{Kov}

Given two ideals $\ical$ and $\jcal $ on ${\omega}$
write $\ical\lerk \jcal$ (see \cite{Ru}) iff there is a function
$f:{\omega}\to {\omega}$ such that
\begin{displaymath}
\ical=\{I\subseteq {\omega}:f^{-1}I\in \jcal\},
\end{displaymath}
and
write $\ical\lerb \jcal$ (see \cite{LaZh}) iff there is a finite-to-one function
$f:{\omega}\to {\omega}$ such that
\begin{displaymath}
\ical=\{I\subseteq {\omega}:f^{-1}I\in \jcal\}.
\end{displaymath}

The following Observations imply that there are $\mc{I}$-MAD families of cardinality $\mf{c}$ for each analytic P-ideal $\mc{I}$.

\begin{obs}\label{ab}
Assume that $\ical$ and  $\jcal$ are ideals on ${\omega}$,
$\ical \lerk \jcal$ witnessed by a function
$f:{\omega}\to {\omega}$.
If  $\acal$ is an  $\ical$-AD family
then $\{f^{-1}A:A\in\acal\}$ is a $\jcal$-AD family.
\end{obs}
\begin{obs}\label{fin_lerb_analp}
$\fini \lerb \ical$ for any analytic $P$-ideal $\ical$.
\end{obs}

\begin{proof}
Let $\ical=\mathrm{Exh}(\varphi)$ for some lower semicontinuous finite submeasure $\varphi$ on $\om$. Since ${\omega}\notin \ical$
we have $\lim_{n\fv\infty}\varphi(\om\bs n)=\varepsilon>0$.
Hence by the lower semicontinuous property of $\varphi$ for each $n>0$ there is $m>n$ such that
$\varphi([n,m))>{\varepsilon}/2$.

So there is a partition  $\{I_n:n<{\omega}\}$ of ${\omega}$
into finite pieces such that $\varphi(I_n)>{\varepsilon}/2$ for each $n\in
{\omega}$. Define the function $f:{\omega}\to{\omega}$
by the stipulation $f''I_n=\{n\}$. Then $f$ witnesses
$\fini \lerb \ical$.
\end{proof}

For any analytic P-ideal $\ical$
denote $\bar{\mf{a}}(\ical)$ the
minimum of the cardinalities of uncountable $\ical$-MAD families.

Clearly $\mf{a}(\ical)>{\omega}$ implies
$\mf{a}(\ical)=\bar{\mf{a}}(\ical)$, especially
$\mf{a}({\ical_h})=\bar{\mf{a}}({\ical_h})$ for summable ideals.

\begin{Thm}\label{a0a}
$\bar{\mf{a}}(\zcal_{\muv})\le\mf{a}$ for each density ideal $\zcal_{\muv}$.
\end{Thm}

\begin{proof}
Let $f:{\omega}\to {\omega}$ be the finite-to-one function defined
by $f^{-1}\{n\}=P_n$ where $\vec{\mu}=\la\mu_n:n\in\om\ra$ and $\mu_n:\mc{P}(P_n)\fv [0,1]$.
Specially $f$ witnesses $\fini \lerb \zcal_\muv$.

Let $\acal$ be an uncountable ($\fini$-)MAD family.
We show that $f^{-1}[\acal]=\{f^{-1}A:A\in \acal\}$ is a $\zcal_\muv$-MAD family.

By Observation \ref{ab}, $f^{-1}[\acal]$ is a $\zcal_\muv$-AD family.

To show the maximality   let $X\in\zcal_\muv^+$ be arbitrary,
$\limsup_{n\fv\infty}{\mu}_n (X\cap P_n)={\varepsilon}>0$. Thus
\begin{displaymath}
J=\{n\in {\omega}:{\mu}_n(X\cap P_n)>{\varepsilon}/2 \}
\end{displaymath}
is infinite. So there is $A\in \acal $ such that
$A\cap J$ is infinite.

Then $f^{-1}A\in f^{-1}[\acal]$ and $X\cap f^{-1}A\in \zcal_\muv^+$
because there are infinitely many $n$ such that we have
$P_n\subseteq f^{-1}A$ and ${\mu}_n(X\cap P_n)>{\varepsilon}/2$.
\end{proof}

\begin{Prob}
{\em Does $\bar{\mf{a}}(\ical)\le\mf{a}$ hold for each analytic P-ideal $\mc{I}$?}
\end{Prob}

\begin{Thm}\label{ba0}
$\mf{b}\le\bar{\mf{a}}(\ical)$ provided that $\ical$ is an
analytic $P$-ideal.
\end{Thm}

\noindent {\em Remark.} If $\mathcal{X}\subset \br {\omega};{\omega};$
is an infinite almost disjoint family then there is a tall ideal $\mathcal{I}$
such that $\mathcal{X}$ is  $\mathcal{I}$-MAD.
So  the Theorem above does not hold for an arbitrary tall
ideal on ${\omega}$.

\begin{proof}
$\ical=\mathrm{Exh}(\varphi)$ for some lower semicontinuous finite submeasure $\varphi$.

Let $\acal$ be an uncountable $\ical$-AD family of cardinality
smaller than $\mf{b}$. We show that $\acal$ is not maximal.

There exists an $\varepsilon>0$ such that the set
\[ \acal_\varepsilon=\big\{A\in\acal:\lim_{n\fv\infty} \varphi(A\bs n)>\varepsilon\big\}\]
is uncountable. Let $\acal'=\{A_n:n\in\om\}\subseteq\acal_\varepsilon$ be a
set of pairwise distinct elements of $\acal_\varepsilon$. We can assume that
these sets are pairwise disjoint.
For each $A\in\acal\bs\acal'$ choose a function $f_A\in\om^\om$ such
that
\begin{itemize}
\item[($*_A$)]
$\varphi\big((A\cap A_n)\setminus f_A(n)\big)< 2^{-n}$
for each $n\in {\omega}$.
\end{itemize}

 Using the assumption $|\acal|<\mf{b}$ there exists a strictly increasing function
$f\in\om^\om$ such that  $f_A\le^*
f$ for each $A\in\acal\bs\acal'$.
For each $n$ pick   $g(n)>f(n)$ such that
$\varphi\big(A_n\cap [f(n), g(n))\big)>\varepsilon$, and
let
\[ X=\mathop{\bigcup}\limits_{n\in\om}\big(A_n\cap [f(n),g(n))\big).\]
Clearly $X\in\zcal_{\muv}^+$ because
for each $n<{\omega}$ there is $m$ such that
$A_m\cap [f(m),g(m))\subs X\bs n$ and so $\varphi(X\setm n)
\ge \varphi\big(A_m\cap [f(m),g(m))\big)> \varepsilon$, i.e. $\lim_{n\fv\infty}\varphi(X\bs n)\fel\varepsilon$.

We have to show that $X\cap A\in\zcal_{\muv}$ for each $A\in\acal$.
If $A=A_n$ for some $n$ then $X\cap A=X\cap A_n=A_n\cap [f(n),g(n))$,
i.e. the intersection is finite.

Assume now that  $A\in\acal\bs\acal'$.
Let ${\delta}>0$. We show that if $k$
is large enough then ${\varphi}((A\cap X)\setm k)<{\delta}$.

There is  $N\in {\omega}$ such that $2^{-N+1}<\de$ and
$f_A(n)\le f(n)$ for each  $n\fel N$.

Let $k$ be so large that $k$ contains the finite set
$\bigcup_{n<N}[f(n),g(n))$.

Now $(X\cap A)\bs k=\bigcup_{n\in\om}\big(A_n\cap A\cap [f(n),g(n))\big)\bs k$ and $\big(A_n\cap A\cap [f(n),g(n))\big)\bs k=\0$ if $n<N$ so
    \begin{multline}\notag
 (X\cap A)\bs k=\bigcup_{n\fel N}\big(A_n\cap A\cap [f(n),g(n))\big)\bs
   k \subseteq\\
\bigcup_{n\fel N}((A_n\cap A)\bs f(n))
\subseteq \bigcup_{n\fel N}((A_n\cap A)\bs f_A(n)).
    \end{multline}
Thus by $(*_A)$ we have
\begin{displaymath}
\varphi((X\cap A)\setm k)\le \sum_{n\fel N}\varphi(A_n\cap A\setm f_A(n))\le
\sum_{n\ge N}\frac{1}{2^n}=2^{-N+1}<{\delta}.
\end{displaymath}

\end{proof}

\section{Cohen-indestructible $\mc{I}$-mad families}\label{sec:cohen}
If $\varphi$ is a lower semicontinuous finite submeasure on $\om$ then clearly $\varphi$ is determined by $\varphi\upharpoonright[\om]^{<\om}$.
Using this observation one can define forcing indestructibility of $\ical$-MAD families for an analytic P-ideal $\ical$.
The following Theorem is a modification of Kunen's proof for existence of Cohen-indestructible MAD family from $\ch$ (see \cite{Ku} Ch. VIII Th. 2.3.). 

\begin{Thm}\label{cohen_indest}
Assume $\ch$. For each analytic P-ideal $\ical$
then there is an uncountable Cohen-indestructible $\ical$-MAD family.
\end{Thm}
\begin{proof}

We will define the uncountable Cohen-indestructible
$\ical$-MAD family $\{A_{\xi}:{\xi}<{\omega}_1\}\subs \ical^+$ by
 recursion on
$\xi\in\om_1$. The family  $\{A_{\xi}:{\xi}<{\omega}_1\}$ 
will be $\mathsf{fin}$-AD as well. Our main concern is that we do have
$\mf{a}(\ical)>{\omega}$ so it is not automatic that
$\{A_{\eta}:{\eta}<{\xi}\}$ is not  maximal for ${\xi}<{\omega}_1$.

Denote $\mbb{C}$ the Cohen forcing.
Let $\mc{I}=\mathrm{Exh}(\varphi)$ be an analytic P-ideal.
Let $\{\la p_\xi,\dot{X}_\xi,\de_\xi\ra:\om\le\xi<\om_1\}$ be an
enumeration
of all triples $\la p,\dot{X},\de\ra$ such that $p\in\mbb{C}$,
$\dot{X}$ is a nice name for a subset of $\om$, and $\de$ is a
positive rational number.

Write ${\varepsilon}=\lim_{n\to \infty}{\varphi}({\omega}\setm n)>0$.
Partition ${\omega}$ into infinite sets $\{A_m:m<{\omega}\}$
such that $\lim_{n\to \infty}{\varphi}(A_m\setm n)={\varepsilon}$
for each $m<{\omega}$.

Assume $\xi\fel\om$ and we have $A_\eta\in\ical^+$ for $\eta<\xi$
such that $\{A_\eta:\eta<\xi\}$ is a $\mathsf{fin}$-AD so especially an $\mc{I}$-AD family.

\smallskip
\noindent{\bf Claim:}  There is $X\in \ical^+$ such that
$|X\cap A_{\zeta}|<{\omega}$ for ${\zeta}<{\xi}$.

\begin{proof}[Proof of the Claim]
Write ${\xi}=\{{\zeta}_i:i<{\omega}\}$.
Recursion on $j\in\om$ we can choose $x_j\in \br A_{\ell_j};<{\omega};$
for some $\ell_j\in {\omega}$
such that
\begin{enumerate}[(i)]
  \item ${\varphi}(x_j)\fel{\varepsilon}/2$,
\item $x_j\cap (\cup_{i\le j}A_{{\zeta}_i})=\emptyset$.
\end{enumerate}
Assume that $\{x_i:i<j\}$ is chosen.
Pick $\ell_j\in {\omega}\setm \{{\zeta}_i:i< j\}$.
Let $m\in {\omega}$ such that
$A_{\ell_j}\cap\cup \{A_{{\zeta}_i}:i\le j\}\subs m$.
Since ${\varphi}(A_{\ell_j}\setm m)\ge {\varepsilon}$ there is
$x_j\in \br A_{\ell_j}\bs m;<{\omega};$ with
${\varphi}(x_j)\ge {\varepsilon}/2$.

Let $X=\cup\{x_j:j<{\omega}\}$. Then $|A_{\zeta}\cap X|<{\omega}$
for ${\zeta}<{\xi}$ and
$\lim_{n\to\infty}(X\setm n)\ge {\varepsilon}/2$.
\end{proof}

If $p_\xi$ does not force (a) and (b) below
then let $A_\xi$  be $X$ from the Claim.
\begin{enumerate}[(a)]
\item $\lim_{n\fv\infty}\check{\varphi}(\dot{X}_\xi\bs
  n)>\check{\de}_\xi$,
\item $\forall$ $\eta<\check{\xi}$
  $\dot{X}_\xi\cap\check{A}_\eta\in\mc{I}$.
\end{enumerate}

Assume $p_\xi\vd$(a)$\wedge$(b).
Let $\{B^\xi_k:k\in\om\}=\{A_\eta:\eta<\xi\}$ and $\{p^\xi_k:k\in\om\}=\{p'\in\mbb{C}:p'\le p_\xi\}$ be enumerations. Clearly for each $k\in\om$ we have
\[ p^\xi_k\vd\lim_{n\fv\infty}\check{\varphi}\big((\dot{X}_{\xi}\bs\cup\{\check{B}^\xi_l:l\le \check{k}\})\bs n\big)>\check{\de}_\xi,\]
so we can choose a $q^\xi_k\le p^\xi_k$ and a finite $a^\xi_k\subseteq
\om$ such that $\varphi(a^\xi_k)>\de_\xi$ and
$q^\xi_k\vd\check{a}^\xi_k\subseteq(\dot{X}_{\xi}\bs\cup\{\check{B}^\xi_l:l\le
\check{k}\})\bs\check{k}$.
Let $A_\xi=\cup\{a^\xi_k:k\in\om\}$. Clearly $A_\xi\in\mc{I}^+$ and $\{A_\eta:\eta\le\xi\}$ is a $\mathsf{fin}$-AD family.

Thus
$\mc{A}=\{A_\xi:\xi<\om_1\}\subseteq\ical^+$ is a $\mathsf{fin}$-AD family.

We show that $\mc{A}$ is a Cohen-indestructible $\mc{I}$-MAD.
Assume otherwise there is a $\xi$ such that
$p_\xi\vd\lim_{n\fv\infty}\check{\varphi}(\dot{X}_\xi\bs
n)>\check{\de}_\xi\wedge\forall$ $\eta<\om_1$
$\dot{X}_\xi\cap\check{A}_\eta\in\mc{I}$,
specially $p_\xi\vd$(a)$\wedge$(b). There is a $p^\xi_k\le p_\xi$ and
an $N$ such that 
$p^\xi_k\vd\check{\varphi}((\dot{X}_\xi\cap\check{A}_\xi)\bs\check{N})<\check{\de}_\xi$.
We can assume $k\fel N$, so
$p^\xi_k\vd\check{\varphi}((\dot{X}_\xi\cap\check{A}_\xi)\bs\check{k})<\check{\de}_\xi$.
By the choice of $q^\xi_k$ and $a^\xi_k$ we have
$q^\xi_k\vd
\check{a}^\xi_k\subseteq(\dot{X}_{\xi}\cap\check{A}_\xi)\bs\check{k}$, so
$q^\xi_k\vd\check{\varphi}((\dot{X}_{\xi}\cap\check{A}_\xi)\bs\check{k})>\check{\de}_\xi$,
contradiction.
\end{proof}

\section{Towers in $\ical^*$}\label{sec:tower}
Let $\ical$ be an ideal on $\om$. A $\subs^*$-decreasing
sequence $\la A_\al:\al<\ka\ra$ is a
{\em tower in $\ical^*$} if (a) it is a tower (i.e. there is no
$X\in \br {\omega};{\omega};$ with $X\subs^* A_{\alpha}$ for
${\alpha}<{\kappa}$), and
(b) $A_\al\in\ical^*$ for $\al<\ka$. Under $\ch$ it is
straightforward to construct towers  in $\mc{I}^*$ for each
tall analytic P-ideal $\mc{I}$.  
The existence of such towers
is consistent with $2^{\omega}>{\omega}_1$ as
well by the  Theorem \ref{tower} below.
Denote 
$\mbb{C}_{\al}$ the standard forcing adding $\al$ Cohen reals by
finite conditions.

\begin{Lem}\label{cohen-tower}
Let $\mc{I}=\mathrm{Exh}(\varphi)$ be a tall analytic P-ideal
in the ground model $V$. Then there is a set $X\in V^{\mbb{C}_1}\cap \mc{I}$
such that $|X\cap S|={\omega}$ for each 
$S\in \br {\omega};{\omega};\cap V$.   
\end{Lem}

\begin{proof}
Since $\mc{I}$ is tall we have
$\lim_{n\fv\infty}\varphi(\{n\})=0$. Fix a
partition $\<I_n:n\in {\omega}\>$ of ${\omega}$
into finite intervals such that
${\varphi}(\{x\})<\frac 1{2^n}$ for $x\in I_{n+1}$ (we can not say
anything
about ${\varphi}(\{x\})$ for $x\in I_0$). Then  $X'\in\mc{I}$
whenever $|X'\cap I_n|\le 1$ for each $n$.

Let  $\{i^n_{k}:k<k_n\}$ be the increasing enumeration of $I_n$.
Our forcing  $\mbb{C}$ adds a
Cohen real $c\in\om^\om$ over $V$. Let
\[ X_\al=\{i^n_k:c(n)\equiv k\;\mathrm{mod}\;{k_n}\}\in
V^{\mbb{C}}\cap\mc{I}.\]
A trivial density argument shows that
 $|X_\al\cap S|=\om$ for each  $S\in V\cap [\om]^\om$.  
\end{proof}

\begin{Thm}\label{tower}
$\vd_{\mbb{C}_{\om_1}}$"There exists a tower in $\ical^*$ for each
  tall analytic P-ideal $\ical$."
\end{Thm}
\begin{proof}
Let $V$ be a countable transitive model and $G$ be a
$\mbb{C}_{\om_1}$-generic filter over $V$.
Let $\mc{I}=\mathrm{Exh}(\varphi)$ be a tall analytic P-ideal in
$V[G]$ with some lower semicontinuous finite submeasure $\varphi$ on
$\om$.
There is a $\de<\om_1$ such that $\varphi\upharpoonright
[\om]^{<\om}\in V[G_\de]$ where $G_\de=G\cap\mbb{C}_\de$, so
we can assume $\varphi\upharpoonright [\om]^{<\om}\in V$.

Work in $V[G]$ recursion on $\om_1$ we construct the tower
$\bar{A}=\la A_\al:\al<\om_1\ra$ in $\mc{I}^*$ such that
$\bar{A}\upharpoonright\al\in V[G_{\al}]$.

Because $\mc{I}$ contains infinite elements we can construct in $V$ a
sequence $\la A_n:n\in\om\ra$ in $\mc{I}^*$ which is strictly
$\subseteq^*$-descending, i.e. $|A_n\bs A_{n+1}|=\om$ for $n\in\om$.
Assume $\la A_\xi:\xi<\al\ra$ are done.

Since $\ical$ is a $P$-ideal there is $A_{\alpha}'\in \ical^*$
with $A_{\alpha}'\subseteq^* A_{\beta}$ for ${\beta}<{\alpha}$.

By lemma \ref{cohen-tower} there is a set $X_{\alpha}\in
V[G_{\al+1}]\cap \mc{I}$ such that $X_{\alpha}\cap S\ne \emptyset$
for each $S\in \br {\omega};{\omega};\cap V[G_{\al}]$.

 Let
$A_\al=A'_\al\bs X_\al\in V[G_{\al+1}]\cap \mc{I}^*$ so
$S\nsubseteq^* A_\al$ for any
$S\in V[G_\al]\cap [\om]^\om$. Hence
$V[G]\models$"$\la A_\al:\al<\om_1\ra$ is a tower in $\mc{I}^*$".
\end{proof}

\begin{Prob}
{\em Do there exist towers in $\mc{I}^*$ for some tall analytic
  P-ideal $\mc{I}$ in $\zfc$?}
\end{Prob}

\section{Unbounding and dominating numbers of ideals}\label{sec:bd}
A {\em supported relation} (see \cite{Vo}) is a triple
$\mc{R}=(A,R,B)$ where $R\subseteq A\times B$, $\dom(R)=A$,
$\ran(R)=B$, and we always assume that for each $b\in B$ there is an
$a\in A$ such that $\la a,b\ra\notin R$.

The {\em unbounding} and {\em dominating numbers of $\mc{R}$}:
\[ \mf{b}(\mc{R})=\min\{|A'|:A'\subseteq A\wedge\forall\ b\in B\
A'\nsubseteq R^{-1}\{b\}\},\]
\[ \mf{d}(\mc{R})=\min\{|B'|:B'\subseteq B\wedge A=R^{-1}B'\}.\]
For example $\mf{b}_\mc{I}=\mf{b}(\om^\om,\le_\mc{I},\om^\om)$ and
$\mf{d}_\mc{I}=\mf{d}(\om^\om,\le_\mc{I},\om^\om)$.
Note that $\mf{b}(\mc{R})$ and $\mf{d}(\mc{R})$ are defined for each
$\mc{R}$, but  in general $\mf{b}(\mc{R})\le\mf{d}(\mc{R})$ does not
hold.

We recall the definition of Galois-Tukey connection of relations.
\begin{definition}(\cite{Vo})
Let $\mc{R}_1=(A_1,R_1,B_1)$ and $\mc{R}_2=(A_2,R_2,B_2)$ be supported relations.
A pair of functions $\phi:A_1\fv A_2$, $\psi:B_2\fv B_1$ is a {\em
  Galois-Tukey connection from $\mc{R}_1$ to $\mc{R}_2$}, in notation
$(\phi,\psi):\mc{R}_1\preceq\mc{R}_2$
if $a_1 R_1\psi(b_2)$ whenever $\phi(a_1)R_2 b_2$.
In a diagram:
\[\begin{CD}
\psi(b_2)\in B_1 @<\psi<< B_2\ni b_2\\
R_1 &\Longleftarrow& R_2\\
a_1\in A_1 @>\phi>> A_2\ni\phi(a_1)
\end{CD}\]
We write $\mc{R}_1\preceq\mc{R}_2$ if there is a Galois-Tukey
connection from $\mc{R}_1$ to $\mc{R}_2$.
If $\mc{R}_1\preceq\mc{R}_2$ and $\mc{R}_2\preceq\mc{R}_1$ also hold
then we say $\mc{R}_1$ and $\mc{R}_2$ are {\em Galois-Tukey
  equivalent}, in notation $\mc{R}_1\equiv\mc{R}_2$.
\end{definition}
\begin{Fact}\label{gt_imp_ineq}
If $\mc{R}_1\preceq\mc{R}_2$
then $\mf{b}(\mc{R}_1)\fel \mf{b}(\mc{R}_2)$ and $\mf{d}(\mc{R}_1)\le
\mf{d}(\mc{R}_2)$.
\end{Fact}

\begin{Thm}
If $\ical \lerb  \jcal$  then
$({\omega}^{\omega},\le_\ical,\om^\om)\equiv
({\omega}^{\omega},\le_\jcal,\om^\om)$.
\end{Thm}

\begin{proof}
Fix a finite-to-one function $f:{\omega}\to {\omega}$ witnessing
$\ical \lerb  \jcal$.

Define $\phi, \psi:{\omega}^{\omega}\to {\omega}^{\omega}$ as follows:

  \begin{gather}\notag
\phi(x)(i)=\max(x''f^{-1}\{i\}),
\\\notag \psi(y)(j)=y(f(j)).
  \end{gather}

\newcommand{\lej}{\le_{\jcal}}
\newcommand{\lei}{\le_{\ical}}

We prove two claims.

  \begin{claim}
$(\phi,\psi):({\omega}^{\omega},\lej,\om^\om)\preceq({\omega}^{\omega},\lei,\om^\om)$.
  \end{claim}

  \begin{proof}[Proof of the claim]
  We show that if $\phi(x)\lei y$ then $x\lej \psi(y)$.
Indeed,   $I=\{i: \phi(x)(i)> y(i)\}\in \ical$.
Assume that $f(j)=i\notin I$.
Then $\phi(x)(i)=\max (x'' f^{-1}\{i\})\le y(i)$.
Since $y(i)=\psi(y)(j)$, so
\begin{displaymath}
x(j)\le \max (x'' f^{-1}\{f(j)\})\le y(f(j))=\psi(y)(j)
\end{displaymath}
Since $f^{-1}I\in \jcal$ this yields $x\lej \psi(y)$.
\end{proof}

\begin{claim}
$(\psi,\phi):({\omega}^{\omega},\lei,\om^\om)\preceq({\omega}^{\omega},\lej,\om^\om)$.
  \end{claim}

  \begin{proof}[Proof of the claim]
We show that if $\psi(y)\lej x$ then $y\lei \phi(x)$.
Assume on the contrary that
$y\not\lei \phi(x)$. Then
$A=\{i\in {\omega}:y(i)>\phi(x)(i)\}\in \ical^+$.
By definition of  $\phi$,
we have $A=\{i:y(i)>\max(x''f^{-1}\{i\})\}$.

Let $B=f^{-1}A\in \jcal^+$.
For  $j\in B$ we have $f(j)\in A$ and so
\begin{displaymath}
\psi(y)(j)=y(f(j))>\phi(x)(f(j))=\max (x'' f^{-1}\{f(j)\})\ge x(j).
\end{displaymath}
 Hence
$\psi(y)\not\lei x$, contradiction.
  \end{proof}
These claims prove the statement of the Theorem, so we are done.
\end{proof}

By Fact \ref{gt_imp_ineq} we have:

\begin{Kov}
If $\ical \lerb  \jcal$ holds then $\mf{b}_\mc{I}=\mf{b}_\jcal$ and
$\mf{d}_\ical=\mf{d}_\jcal$.
\end{Kov}

By Observation \ref{fin_lerb_analp} this yields:

\begin{Kov}\label{cor:equiv}
If $\ical$ is an analytic $P$-ideal  then
$({\omega}^{\omega},\le^*,\om^\om)\equiv
({\omega}^{\omega},\le_\jcal,\om^\om)$,
and
$\mf{b}_\ical=\mf{b}$ and
$\mf{d}_\ical=\mf{d}$.
\end{Kov}

\section{$\ical$-bounding and $\ical$-dominating forcing notions}\label{sec:forcing}
\begin{definition}
Let $\ical$ be a Borel ideal on $\om$. A forcing notion $\PP$ is
{\em $\ical$-bounding} if
\[ \vd_\PP\forall\; A\in \ical\;\exists\;B\in\ical\cap V\;A\subseteq
B;\] $\PP$ is {\em $\ical$-dominating} if
\[ \vd_\PP\exists\; B\in \ical\;\forall\;A\in\ical\cap V\;A\subseteq^*
B.\]
\end{definition}
\begin{Thm}\label{ibounddom}
Let $\ical$ be a tall analytic $P$-ideal.
 If $\PP$ is
$\ical$-bounding then $\PP$ is $\om^\om$-bounding as well; if $\PP$
is $\ical$-dominating then $\PP$ adds dominating reals.
\end{Thm}
\begin{proof}
Assume that $\ical=\mathrm{Exh}({\varphi})$ for some lower semicontinuous finite
submeasure ${\varphi}$.
For $A\in \ical$ let
\begin{displaymath}
  d_A(n)=\min\bigl\{k\in {\omega}:{\varphi}(A\setm k)<{2^{-n}}\bigr\}.
\end{displaymath}
Clearly if $A\subseteq B\in\mc{I}$ then $d_A\le d_B$.

It is enough to show that $\{d_A:A\in\mc{I}\}$ is cofinal in
$\la\om^\om,\le^*\ra$. Let $f\in\om^\om$.
Since $\ical$ is a tall ideal  we have $\lim_{k\to\infty}{\varphi}(\{k\})=0$
but $\lim_{m\to\infty}({\omega}\setm m)=\varepsilon>0$.
Thus for all but finite $n\in\om$ we can choose
a finite set $A_n\subseteq\om\bs f(n)$ such that $2^{-n}\le {\varphi}(A_n)<
2^{-n+1}$ so $A=\cup\{A_n:n\in\om\}\in\ical$ and $f\le^* d_A$.

Why? We can assume if $k\fel f(n)$ then $\varphi(\{k\})<2^{-n}$.
Let $n$ be so large such that $2^{-n}<\varepsilon$.
Now if there is no a suitable $A_n$ then $\varphi(\om\bs f(n))\le 2^{-n}<\eps$, contradiction.
\end{proof}

The converse of the first implication of Theorem \ref{ibounddom} is
not true by the following Proposition.

\begin{Prop}\label{random}
The random forcing is not $\ical$-bounding for any tall summable and tall density ideal $\mc{I}$.
\end{Prop}

\begin{proof}

Denote $\mbb{B}$ the random forcing and $\lam$ the Lebesgue-measure.

If 
$\mc{I}=\mc{I}_h$ is a tall summable ideal then 
we can chose pairwise disjoint sets $H(n)\in [\om]^\om$ such that
$\sum_{l\in H(n)}h(l)=1$ and $\max\{h(l):l\in H(n)\}<2^{-n }$ for each
$n\in\om$.
Let $H(n)=\{l^n_k:k\in\om\}$. For each $n$ fix a partition
$\{[B^n_k]:k\in\om\}$ of $\mbb{B}$ such that $\lam(B^n_k)=h(l^n_k)$
for each $k\in\om$.
Let $\dot{X}$ be a $\mbb{B}$-name such that
$\vd_\mbb{B}\dot{X}=\{\check{l}^n_k:\check{[B^n_k]}\in\dot{G}\}$.
Clearly $\vd_{\mbb{B}}\dot{X}\in\mc{I}_h$. $\dot{X}$ shows that
$\mbb{B}$ is not $\mc{I}_h$-bounding.

Assume on the contrary that there is a $[B]\in\mbb{B}$ and an
$A\in\mc{I}_h$ such that $[B]\vd\dot{X}\subseteq\check{A}$.
There is an $n\in\om$ such that
\[ \mathop{\sum}\limits_{l^n_k\in A}\lam(B^n_k)=\sum_{l^n_k\in A}h(l^n_k)<\lam(B).\]
Choose a $k$ such that $l^n_k\notin A$ and $[B^n_k]\wedge [B]\ne
[\0]$. We have $[B^n_k]\wedge [B]\vd
\check{l}^n_k\in\dot{X}\bs\check{A}$, contradiction.\\

If $\mc{I}=\zcal_{\vec{\mu}}$ is a tall density ideal then 
for each $n$ fix a partition $\{[B^n_k]:k\in P_n\}$ of $\mbb{B}$ such
that $\lam(B^n_k)=\mu_n(\{k\})$ for each $k$.
Let $\dot{X}$ be a $\mbb{B}$-name such that
$\vd_\mbb{B}\dot{X}=\{\check{k}:\check{[B^n_k]}\in\dot{G}\}$. Clearly
$\vd_{\mbb{B}}\dot{X}\in\mc{Z}_{\vec{\mu}}$. $\dot{X}$ shows that
$\mbb{B}$ is not $\mc{Z}_{\vec{\mu}}$-bounding.

Assume on the contrary that there is a $[B]\in\mbb{B}$ and an
$A\in\mc{Z}_{\vec{\mu}}$ such that $[B]\vd\dot{X}\subseteq\check{A}$.
There is an $n\in\om$ such that
\[ \sum_{k\in A\cap P_n}\lam(B^n_k)=\mu_n(A\cap P_n)<\lam(B).\]
Choose a $k\in P_n\bs A$ such that $[B^n_k]\wedge [B]\ne [\0]$. We
have $[B^n_k]\wedge [B]\vd \check{k}\in\dot{X}\bs\check{A}$,
contradiction.
\end{proof}
The converse of the second implication of Theorem \ref{ibounddom} is
not true as well: the Hechler forcing is a counterexample according to
the following 
Theorem.
\begin{Thm}\label{sigma}\label{tm:sigma}
If $\PP$ is $\sigma$-centered then $\PP$ is not $\ical$-dominating for
any tall analytic $P$-ideal  $\ical$.
\end{Thm}

\begin{proof}
Assume that $\ical=\mathrm{Exh}({\varphi})$ for some lower semicontinuous
finite submeasure ${\varphi}$.
Let ${\varepsilon}=\lim_{n\to \infty}{\varphi}({\omega}\setm n)>0$.

Let $\PP=\cup\{C_n:n\in\om\}$ where $C_n$ is centered for each $n$.
Assume on the contrary that $\vd_\PP\dot{X}\in\mc{I}\wedge\forall$
$A\in\mc{I}\cap V$ $A\subseteq^* \dot{X}$ for some $\PP$-name
$\dot{X}$.

 For each $A\in\ical$
choose a $p_A\in\PP$ and a $k_A\in\om$ such that
\begin{equation}\tag{$\circ$}
p_A\vd\check{A}\bs\check{k}_A\subseteq\dot{X}
\wedge{\varphi}(\dot{X}\setm\check{k}_A)<{\varepsilon}/2.
\end{equation}

For each $n,k\in\om$ let $\mc{C}_{n,k}=\{A\in\ical:p_A\in C_n\wedge k_A=k\}$,
and let $B_{n,k}=\bigcup\mc{C}_{n,k}$.
We show that for each $n$ and
$k$
$${\varphi}(B_{n,k}\setm k)\le {\varepsilon}/2.$$

Assume indirectly ${\varphi}(B_{n,k}\bs k)>{\varepsilon}/2$
for some $n$ and $k$.
There is a $k'$ such that ${\varphi}( B_{n,k}\cap [k,k'))>{\varepsilon}/2$ and
  there is a finite $\mc{D}\subseteq\mc{C}_{n,k}$ such that
  $B_{n,k}\cap [k,k')=(\cup\mc{D})\cap [k,k')$. Choose a common
      extension
$q$ of $ \{p_A:A\in\mc{D}\}$.
Now we have
$q\vd\cup\{A\bs\check{k}:A\in\check{\mc{D}}\}\subseteq\dot{X}$ and so
\[ q\vd {\varepsilon}/2<{\varphi}( \check{B}_{n,k}\cap
  [\check{k},\check{k}'))=
{\varphi}((\cup\check{\mc{D}})\cap [\check{k},\check{k}'))
\le {\varphi}(\dot{X}\cap [\check{k},\check{k}'))
\le {\varphi}(\dot{X}\bs\check{k}),\]
which contradicts $(\circ)$.

So for each $n$ and $k$ the set ${\omega}\setm B_{n,k}$ is infinite, so ${\omega}\setm B_{n,k}$
contains an infinite $D_{n,k}\in \ical$.
Let $D\in \ical$ such that $D_{n,k}\subs^* D$ for each $n,k\in {\omega}$.

Then there is no $n,k$ such that $D\subs^* B_{n,k}$. Contradiction.
\end{proof}
By this Theorem an by Lemma \ref{cohen-tower} 
the Cohen forcing is neither $\ical$-dominating 
nor $\ical$-bounding
for any tall analytic P-ideal $\ical$.

Finally in the rest of the paper we compare the Sacks property and the
$\mc{I}$-bounding property.
\begin{Thm}\label{sacksibound}
If $\PP$ has the Sacks property then $\PP$ is $\mc{I}$-bounding for
each analytic P-ideal $\mc{I}$.
\end{Thm}

\begin{proof}
Let $\mc{I}=\mathrm{Exh}(\varphi)$.
Assume $\vd_\PP\dot{X}\in\mc{I}$. Let $d_{\dot{X}}$ be a $\PP$-name
for an element of $\om^\om$
such that $\vd_\PP
d_{\dot{X}}(\check{n})=\min\{k\in\om:\varphi(\dot{X}\bs
k)<2^{-\check{n}}\}$.
We know that $\PP$ is $\om^\om$-bounding. If $p\vd
d_{\dot{X}}\le\check{f}$ for some strictly increasing $f\in\om^\om$
then by the Sacks property there is a $q\le p$ and a slalom
$S:\om\fv\big[[\om]^{<\om}\big]^{<\om}$, $|S(n)|\le n$ such that
\[q\vd\forall^{\infty}\,n\;\dot{X}\cap [f(n),f(n+1))\in S(n).\]
Now let
\[ A=\bigcup_{n\in\om}\{D\in S(n):\varphi(D)<2^{-n}\}.\]
$A\in\mc{I}$ because $\varphi(A\bs f(n))\le\sum_{k\fel n}\varphi(A\cap [f(k),f(k+1))\le\sum_{k\fel n}\frac{k}{2^k}$. Clearly $q\vd\dot{X}\subseteq^*\check{A}$.
\end{proof}

A supported relation $\mc{R}=(A,R,B)$ is called {\em
Borel-relation} iff there is a Polish space $X$ such that $A,
B\subseteq X$ and $R\subseteq X^2$ are Borel sets. Similarly a
Galois-Tukey connection $(\phi,\psi):\mc{R}_1\preceq\mc{R}_2$
between Borel-relations is called {\em Borel GT-connection} iff
$\phi$ and $\psi$ are Borel functions. To be Borel-relation and
Borel GT-connection is absolute for transitive models containing all
relevant codes.

Some important Borel-relation:

(A): $(\mc{I},\subseteq,\mc{I})$ and $(\mc{I},\subseteq^*,\mc{I})$
for a Borel ideal $\mc{I}$.

(B): Denote $\mathrm{Slm}$ the set of {\em slaloms on $\om$}, i.e.
$S\in\mathrm{Slm}$ iff $S:\om\fv [\om]^{<\om}$ and $|S(n)|=2^n$ for
each $n$. Let $\sqsubseteq$ and $\sqsubseteq^*$ be the following
relations on $\om^\om\times\mathrm{Slm}$:
\[ f\sqsubseteq^{(*)} S\iff\forall^{(\infty)}\,n\in\om\,f(n)\in S(n).\]
The supported relations $(\om^\om,\sqsubseteq,\mathrm{Slm})$ and $(\om^\om,\sqsubseteq^*,\mathrm{Slm})$ are  Borel-relations.

(C): Denote $\ell^+_1$ the set of positive summable series. Let
$\le$ be the coordinate-wise and $\le^*$ the almost everywhere
coordinate-wise ordering on $\ell^+_1$. $(\ell^+_1,\le,\ell^+_1)$
and $(\ell^+_1,\le^*,\ell^+_1)$ are Borel-relations.

\begin{definition}
Let $\mc{R}=(A,R,B)$ be a Borel-relation. A forcing notion $\PP$ is
{\em $\mc{R}$-bounding} if
\[ \vd_\PP\forall\,a\in A\,\exists\,b\in B\cap V\, aRb;\]
{\em $\mc{R}$-dominating} if
\[ \vd_\PP\exists\,b\in B\,\forall\,a\in A\cap V\,aRb.\]
\end{definition}

For example the property $\mc{I}$-bounding/dominating is the same as $(\mc{I},\subseteq^{*},\mc{I})$-bounding/dominating.

We can reformulate some classical properties of forcing notions:
\begin{align*}
\om^\om\text{-bounding} && \equiv &&(\om^\om,\le^{(*)},\om^\om)\text{-bounding}\\
\text{adding dominating reals} && \equiv && (\om^\om,\le^*,\om^\om)\text{-dominating}\\
\text{Sacks property} && \equiv && (\om^\om,\sqsubseteq^{(*)},\mathrm{Slm})\text{-bounding}\\
\text{adding a
slalom capturing} && \equiv &&
(\om^\om,\sqsubseteq^*,\mathrm{Slm})\text{-dominating}\\
\text{all ground model reals} && &&
\end{align*}

If $\mc{R}=(A,R,B)$ is a supported relation then let $\mc{R}^\perp=(B,\neg R^{-1},A)$ where $b(\neg R^{-1}) a$ iff not $a Rb$. Clearly $(\mc{R}^\perp)^\perp=\mc{R}$ and $\mf{b}(\mc{R})=\mf{d}(\mc{R}^\perp)$. Now if $\mc{R}$ is a Borel-relation then $\mc{R}^\perp$ is a Borel-relation too, and a forcing notion is $\mc{R}$-bounding iff it is not $\mc{R}^\perp$-dominating.

\begin{Fact}\label{force_gt_imp}
Assume $\mc{R}_1\preceq\mc{R}_2$ are Borel-relations with Borel GT-connection and $\PP$ is a forcing notion. If $\PP$ is $\mc{R}_2$-bounding/dominating then $\PP$ is $\mc{R}_1$-bounding/dominating.
\end{Fact}

By Corollary \ref{cor:equiv} this yields
\begin{Kov}\label{cor:dom}
For each analytic $P$-ideal $\mc{I}$
(1) a poset $\PP$  is $\le_{\mc{I}}$-bounding iff it is ${\omega}^{\omega}$-bounding,
(2) forcing with a poset $\PP$ adds $\le_{\mc{I}}$-dominating reals iff
this forcing adds
dominating reals.
\end{Kov}

We will use the following Theorem.
\begin{Thm}{\em (\cite{Fr} 526B, 524I)}\label{ketFremlin} There are Borel GT-connections
$(\mc{Z},\subseteq,\mc{Z})\preceq (\ell^+_1,\le,\ell^+_1)$ and
$(\ell^+_1,\le^*,\ell^+_1)\equiv (\om^\om,\sqsubseteq^*,\mathrm{Slm})$.
\end{Thm}
Note that there is no any Galois-Tukey connection from $(\ell^+_1,\le,\ell^+_1)$ to $(\mc{Z},\subseteq,\mc{Z})$ so they are not GT-equivalent (see \cite{LoVe}) Th. 7.).
\begin{Kov}
If $\PP$ adds a slalom capturing all ground model reals then $\PP$ is $\mc{Z}$-dominating.
\end{Kov}
\begin{proof}
By Fact \ref{force_gt_imp} and Theorem \ref{ketFremlin} adding slalom is the same as $(\ell^+_1,\le^*,\ell^+_1)$-dominating. Let $\dot{x}$ be a $\PP$-name such that $\vd_\PP\dot{x}\in\ell^+_1\wedge\forall$ $y\in\ell^+_1\cap V$ $y\le^*\dot{x}$. Moreover let $\dot{X}$ be a $\PP$-name such that $\vd_\PP\dot{X}=\{z\in\ell^+_1:|z\bs\dot{x}|<\om$, $\forall$ $n$ $(z(n)\ne \dot{x}(n)\Ra z(n)\in\om)\}$.
Let $(\phi,\psi): (\mc{Z},\subseteq,\mc{Z})\preceq (\ell^+_1,\le,\ell^+_1)$ be a Borel GT-connection.
Now if $\dot{A}$ is a $\PP$-name such that $\vd_\PP\forall$ $z\in\dot{X}$ $\psi(z)\subseteq^*\dot{A}$ then $\dot{A}$ shows that $\PP$ is $\mc{Z}$-dominating.
\end{proof}

Denote $\mbb{D}$ the dominating forcing and $\mbb{LOC}$ the Localization forcing.

\begin{obs}
If $\mc{I}$ is an arbitrary analytic P-ideal then  two step iteration
$\mathbb{D}*\mbb{LOC}$ is $\mc{I}$-dominating.
\end{obs}

Indeed, let $\mc{I}\in V\subseteq M\subseteq N$ be transitive models,
$d\in M\cap\om^\om$ be strictly increasing and dominating over $V$,
and $S\in N$, $S:\om\fv\big[[\om]^{<\om}\big]^{<\om}$, $|S(n)|\le n$ a
slalom which captures all reals from $M$. Now if
\[ X_n=\cup\{A\in S(n)\cap \mc{P}([d(n),d(n+1)):\varphi(A)<2^{-n}\}\]
then it is easy to see that $Y\subseteq^*\cup\{X_n:n\in\om\}\in
\ical\cap N$ for each $Y\in V\cap \mc{I}$.

\begin{Prob}
{\em For which analytic P-ideal $\mc{I}$ does
$(\mc{I},\subseteq^{(*)},\mc{I})\preceq
(\ell^+_1,\le^{(*)},\ell^+_1)$ hold, or ``adding slaloms'' imply $\mc{I}$-dominating, or at least $\mbb{LOC}$ is $\mc{I}$-dominating?}
\end{Prob}

\begin{Prob}
{\em Does $\zcal$-dominating (or $\mc{I}$-dominating) imply adding slaloms?}
\end{Prob}

We will use the following deep result of Fremlin to prove Theorem \ref{zb-s}.
\begin{Thm}\label{balszorzat}{\em (\cite{Fr} 526G)}
There is a family $\{P_f:f\in\om^\om\}$ of Borel subsets of $\ell^+_1$ such that the following hold:
\begin{itemize}
\item[(i)] $\ell^+_1=\cup\{P_f:f\in\om^\om\}$,
\item[(ii)] if $f\le g$ then $P_f\subseteq P_g$,
\item[(iii)] $(P_f,\le,\ell^+_1)\preceq(\mc{Z},\subseteq,\mc{Z})$ with a Borel GT-connection for each $f$.
\end{itemize}
\end{Thm}

\begin{Thm}\label{zb-s}
$\PP$ is $\mc{Z}$-bounding iff $\PP$ has the Sacks property.
\end{Thm}

\begin{proof}
Let $\{P_f:f\in\om^\om\}$ be a family satisfying (i), (ii), and (iii) in Theorem \ref{balszorzat}, and fix Borel GT-connections $(\phi_f,\psi_f):(P_f,\le,\ell^+_1)\preceq (\mc{Z},\subseteq,\mc{Z})$ for each $f\in\om^\om$.
Assume $\PP$ is $\mc{Z}$-bounding and $\vd_\PP\dot{x}\in\ell^+_1$. $\PP$ is $\om^\om$-bounding by Theorem \ref{ibounddom} so using (ii) we have $\vd_\PP\ell^+_1=\cup\{P_f:f\in\om^\om\cap V\}$. We can choose a $\PP$-name $\dot{f}$ for an element of $\om^\om\cap V$ such that
$\vd_\PP\dot{x}\in P_{\dot{f}}$. By $\mc{Z}$-bounding property of $\PP$ there is a $\PP$-name $\dot{A}$ for an element of $\mc{Z}\cap V$
such that $\vd_\PP\phi_{\dot{f}}(\dot{x})\subseteq\dot{A}$, so $\vd_\PP\dot{x}\le\psi_{\dot{f}}(\dot{A})\in\ell^+_1\cap V$. So we have $\PP$ is $(\ell^+_1,\le^{(*)},\ell^+_1)$-bounding. By Theorem \ref{ketFremlin} and Fact \ref{force_gt_imp} $\PP$ has the Sacks property.

The converse implication was proved in Theorem \ref{sacksibound}.

\end{proof}

\begin{Prob}
{\em Does the $\mc{I}$-bounding property imply the Sacks property for each tall analytic P-ideal $\mc{I}$?}
\end{Prob}

\end{document}